\documentclass[a4wide]{article}

\usepackage[dvips]{epsfig} 
\usepackage[dvips]{hyperref}
\usepackage[latin1]{inputenc}
\usepackage{amsmath}
\usepackage{amsfonts}
\usepackage{amssymb}

\title{Iterative schemes for computing fixed points of nonexpansive mappings in Banach spaces}
\author{Jean-Philippe Chancelier \thanks{Cermics, \'Ecole Nationale des Ponts et Chauss\'ees, 6 et 8 avenue Blaise Pascal, 77455, Marne la Vall\'ee, Cedex, France}}

\newtheorem{thm}{Theorem}
\newtheorem{cor}[thm]{Corollary}
\newtheorem{lem}[thm]{Lemma}

\newtheorem{rem}[thm]{Remark}

\begin{document}

\maketitle

\begin{abstract}
Let $X$ be a real Banach space with a normalized duality mapping uniformly norm-to-weak$^\star$ 
continuous on bounded sets or a reflexive Banach space which admits a weakly continuous duality 
mapping $J_{\Phi}$ with gauge $\phi$. Let $f$ be an {\em $\alpha$-contraction} and $\{T_n\}$ a 
sequence of nonexpansive mapping, we study the strong convergence of explicit iterative schemes 
\begin{equation}
  x_{n+1} = \alpha_n f(x_n) + (1-\alpha_n) T_n x_n 
\end{equation}
with a general theorem and then recover and improve some specific cases studied in the literature 
\cite{xu-2004,xu-2005,song-chen,song-chen-1,chen-zhang-fan,kimura}. 
\end{abstract}



\def\Fix{\mathop{\normalfont Fix}}
\def\defpar{\stackrel{\mbox{\tiny def}}{=}}
\def\argmax{\mathop{\mbox{\rm Argmax}}}
\def\argmin{\mathop{\mbox{\rm Argmin}}}
\def\bbP{{\mathbb P}} 
\def\bbC{{\mathbb C}} 
\def\bbE{{\mathbb E}} 
\def\E{{\cal E}}
\def\F{{\cal F}}
\def\H{{\cal H}}
\def\V{{\cal V}}
\def\W{{\cal W}}
\def\U{{\cal U}}
\def\R{{\mathbb R}}
\def\RB{{\mathbb R}}
\def\N{{\mathbb N}}
\def\M{{\mathbb M}}
\def\S{{\mathbb S}}
\def\bbR{{\mathbb R}} 
\def\bbI{{\mathbb I}} 
\def\bbU{{\mathbb U}} 
\def\Tad{{{\cal T}_{\mbox{\tiny ad}}}}
\def\texte#1{\quad\mbox{#1}\quad}
\def\Proba#1{\bbP\left\{ #1 \right\}} 
\def\Probax#1#2{{\bbP}_{#1}\left\{ #2 \right\}} 
\def\ProbaU#1#2{{\bbP}^{#1} \left\{ #2 \right\}} 
\def\ProbaxU#1#2#3{{\bbP}^{#1}_{#2} \left\{ #3 \right\}} 
\def\valmoy#1{\bbE\left[ #1 \right]}
\def\valmoyDebut#1{\bbE [ #1 } 
\def\valmoyFin#1{ #1 ]} 
\def\valmoyp#1#2{\bbE_{#1}\left[ #2 \right]}
\def\valmoypDebut#1#2{\bbE_{#1} \left[ #2 \right.} 
\def\valmoypFin#1{ \left. #1 \right]} 
\def\valmoypU#1#2#3{\bbE_{#1}^{#2}\left[ #3 \right]}
\def\norminf#1{ {\Vert #1 \Vert}_{\infty}}
\def\norm#1{ {\Vert #1 \Vert}}
\def\Hun{${\text{\bf H}}_1$}
\def\Hdeux{${\text{\bf H}}_2$}
\def\Htrois{${\text{\bf H}}_3$}
\def\psca#1{\left< #1 \right>}
\def\slim{\sigma\mbox{-}\lim}

\newenvironment{myproof}{{\small{\it Proof~:}}}{\hfill$\Box$\normalsize
\\\smallskip}

\section{Introduction and preliminaries} 

Let $X$ be a real Banach space, $C$ a nonempty closed convex subset 
of $X$. Recall that a mapping $T : C \mapsto  C$ is {\em nonexpansive} if
$\norm{T(x) - T(y) } \le \norm{x -y}$ for all $x$, $y\in C$ and a mapping 
$f : C \mapsto  C$ is an {\em $\alpha$-contraction} if there exists 
$\alpha\in(0,1)$ such that $\norm{f(x) - f(y) } \le \alpha \norm{x -y}$  
for all $x$, $y\in C$. 

We denote by $\Fix(T)$ the set of fixed points of $T$, that is 
\begin{equation}
  \Fix(T) \defpar  \left\{x \in C\,  :\,  T x = x \right\}
\end{equation}
and $\Pi_C$ will denote the collection of contractions on $C$. 

Let $X$ be a real Banach space. The (normalized) duality map $J : X \mapsto X^\star$, where $X^\star$
 is the dual space of $X$, is defined by~:
\begin{equation*}
  J(x) \defpar \left\{ x^\star \in X^\star\, :\, \psca{x, x^\star} = \norm{x}^2 = \norm{x^\star}^2 \right\}
\end{equation*}
and there holds the inequality
\begin{equation*}
\norm{x + y}^2 \le \norm{x}^2 + 2 \psca{y , j (x + y)} \text{ where } x,y\in X \text{ and } 
j (x + y) \in J(x + y)\,. 
\end{equation*}

Recall that if $C$ and $F$ are nonempty subsets of a Banach space $X$ such that $C$ is nonempty
closed convex and $F \subset C$, then a map $R : C \mapsto  F$  is called a 
{\em retraction} from $C$ onto
$F$ if $R(x) = x$  for all $x \in F$. A retraction $R : C \mapsto F$ is {\em sunny} provided
$R(x + t (x - R(x))) = R(x)$ for all $ x \in C$ and $ t\ge 0$ whenever $x + t (x - R(x)) \in C$.
A {\em sunny nonexpansive retraction} is a sunny retraction, which is also 
nonexpansive.

Suppose that $F$ is the non empty fixed point set of a nonexpansive mapping $T : C \mapsto C$,
that is $F=\Fix{T}\ne \emptyset$ and assume that $F$ is closed. 
For a given $u\in C$ and every $t \in (0, 1)$ there exists a fixed 
point, denoted $x_t$, of the $(1-t)$-contraction $t u + (1-t)T$. 
Then we define $Q : C \mapsto F=\Fix(T)$ by $Q(u) \defpar \slim_{t \to 0} x_t$ when 
this limit exists ($\slim$ denotes the strong limit). $Q$ will also be denoted by 
$Q_{\Fix(T)}$ when necessary and note that it is easy to check that, when it exists, 
$Q$ is a nonexpansive retraction. 

Consider now $f$ an {\em $\alpha$-contraction}, then $Q_{\Fix(T)} \circ f$ is also an 
{\em $\alpha$-contraction} and admits therefore a unique fixed point 
$\tilde{x} = Q_T \circ f(\tilde{x})$. We define by ${\mathbf Q}(f)$ or 
${\mathbf Q}_{\Fix(T)}(f)$ the mapping ${\mathbf Q}(f) : \Pi_C \rightarrow \Fix{(T)}$ 
such that~: 
\begin{equation}
  {\mathbf Q}(f) \defpar \tilde{x} \quad \text{  where } \quad\tilde{x} = (Q_{\Fix(T)} \circ f) (\tilde{x}). 
  \label{bfq}
\end{equation}

For $t\in(0,1)$ we can also find a fixed point, denoted $x_t^f$ of the $(1-(1-t)\alpha)$-contraction 
$t f + (1-t)T$ and if $\lim_{t \to 0} x_t^f$ is well defined we can define a mapping 
$\widetilde{\mathbf Q}: \Pi_C \mapsto \Fix(T)$ by~:

\begin{equation}
  \widetilde{\mathbf Q}(f) \defpar \lim_{t \to 0} x_t^f \quad \text{  where } \quad 
   x_t^f = t f( x_t^f) + (1-t)T x_t^f
  \label{tbfq}
\end{equation}

We then gather know theorems under which $Q$, ${\mathbf Q}$ and $\widetilde{\mathbf Q}$ are defined 
and give relations between them. 

When $X$ is a uniformly smooth Banach space, denoted by ${\cal B}_{\mbox{us}}$, It is known 
\cite[Theorem 4.1]{xu-2004} that $\widetilde{\mathbf Q}(f)$ is well defined and equal to 
${\mathbf Q}(f)$ and $\tilde{x}={\mathbf Q}(f)$ is characterized by~: 
\begin{equation} 
  \psca{\tilde{x} - f(\tilde{x}), J( \tilde{x} -p)} \le 0 \text{ for all } p\in F=\Fix{(T)}.
  \label{viequation}
\end{equation}
A special case is when $f$ is a constant function ${\mathbf u}(x)=u$. Then \cite[Theorem 4.1]{xu-2004} 
shows that $Q$ is well defined and that $Q(u) = {\mathbf Q}({\mathbf u})= P_{\Fix{T}} u$ 
(where $P_S$ is the metric projection on $S$). If $X$ is a smooth Banach space, 
$R : C \mapsto F$ is a sunny nonexpansive retraction \cite{goebel-reich} if and only if the following 
inequality holds~: 
\begin{equation}
\psca{x - R x, J (y - R x)} \le 0 \text{ for all } x\in C \text{ and } y \in F.
\label{sunnyretract}
\end{equation}
$Q$ is thus the unique sunny non expansive retraction from $C$ to $\Fix{T}$. \cite[Theorem 4.1]{xu-2004}  
was already known in the case f constant and in the context of Hilbert spaces 
\cite[Theorem 3.1]{xu-2004} and \cite[Theorem 2.1]{moudafi}. 

The same existence and characterization results can be found firstly
when $X$ is a reflexive Banach space which admits a weakly continuous duality mapping $J_{\Phi}$ with gauge $\phi$, 
denoted by ${\cal B}_{\text{rwsc}}$, in \cite[Theorem 3.1]{xu-2006} (with $f$ constant) and \cite[Theorem 2.2]{song-chen-1} 
(where $J$ is the (normalized) duality mapping). Note that the limitation of $f$ constant in \cite{xu-2006} can be relaxed 
with \cite{suzuki}. 
Secondly when $X$ is a reflexive and a strictly convex Banach space with a uniformly Gâteaux differentiable norm,  
denoted by ${\cal B}_{\text{rug}}$, \cite[Theorem 3.1]{song-chen}. 
Note that in this three Banach spaces cases listed here the normalized duality mapping is shown to be single valued. 

The aim of this paper is to study the strong convergence of iterative schemes~: 
\begin{equation}
  x_{n+1} = \alpha_n f(x_n) + (1-\alpha_n) T_n x_n 
  \label{main}
\end{equation}
when $X$ can be a ${\cal B}_{\mbox{us}}$, or a ${\cal B}_{\text{rwsc}}$, or a ${\cal B}_{\text{rug}}$ real 
Banach space and $\{T_n\}$ is a sequence of nonexpansive mappings which share at least a common fixed point. 
We give a general framework to show that $\{x_n\}$ will converge strongly to $\tilde{x}$ 
where $\tilde{x}$ is the unique solution of \eqref{viequation} for a fixed nonexpansive mapping 
$T$ related to the sequence $\{T_n\}$. 
The key ingredient is the fact that Lemma \ref{limsuplemma} given in section \ref{secbib} is valid 
in the three previous context. Then we show that by specifying the sequence $T_n$ we can recover 
and extend some known convergence theorems 
\cite{xu-2004,xu-2005,song-chen,song-chen-1,chen-zhang-fan,kimura}. 
Note also that in equation \eqref{main}, $f$ is an $\alpha$-contraction, 
but following \cite{suzuki} it is easy to show that $f$ can be replaced by a Meir-Keeler 
contraction (Lemma \ref{mkclem} in section \ref{secbib} is devoted to this extension).
The paper is organized as follows~: a main theorem is proved in section \ref{secbib} using a set of lemmas 
which are postponed to the last section of the paper and which are verbatim or slight extensions of know 
results. Then in a collection of subsections, known convergence theorems are revisited with shorter proofs. 

\section{Main theorem}

In the sequel a {\em ${\cal B}$ real Banach space}, will denote when not specifically stated a 
real Banach space with a normalized duality mapping uniformly norm-to-weak$^\star$ continuous on bounded sets 
(which is the case for ${\cal B}_{\mbox{us}}$ or ${\cal B}_{\text{rug}}$) or a reflexive Banach space which admits 
a weakly continuous duality mapping $J_{\Phi}$ with gauge $\phi$ (${\cal B}_{\text{rwsc}}$). 

\label{mainsec}
\def\hypothd#1#2{$\mbox{\bf H}_{\mathbf{#1,#2}}$}
\def\hypothdp#1#2{$\mbox{\bf H}^{'}_{\mathbf{#1,#2}}$}
\def\hypothu#1{$\mbox{\bf H}_{\mathbf{#1}}$}

\hypothd{1}{N}: For a fixed given $N \ge 1$ and a given sequence $\{\alpha_n\}$, 
a sequence of mappings $\{T_n\}$ will be said to verify \hypothd{1}{N}, if for a 
given bounded sequence $\{z_n\}$, we have 
\begin{equation}
  \norm{ (1-\alpha_{n+N}) T_{n+N} z_{n} - (1 - \alpha_{n}) T_n z_n} \le \delta_n M 
  \label{hunequation}
\end{equation}
with either $(i)$ $\sum_{0}^\infty |\delta_n| < \infty$ or 
$(i')$ $\limsup_{n\to \infty} \delta_n / \alpha_{n} \le 0$ 
and $M$ a constant.

\begin{rem}Note that using Lemma \ref{lemmaconvjpc} $\{\delta_n\}$ can be replaced by $\{\mu_n + \rho_n\}$ 
where $\{\mu_n\}$ satisfies $(i)$ and $\{\rho_n\}$ satisfies $(i')$.
\end{rem}

\begin{rem}Note that when $\alpha_n \in (0,1)$ we have~: 
\label{hypoth1rem} 
\begin{equation}
  \norm{ (1-\alpha_{n+N}) T_{n+N} z_{n} - (1 - \alpha_{n}) T_n z_n} 
  \le | \alpha_{n+N} - \alpha_n| \norm{T_{n+N} z_n} + \norm{ T_{n+N} z_n - T_n z_n}.
  \label{hunequationrem}
\end{equation}
Thus, when $\{\alpha_n\}$ satisfies \hypothd{3}{N} (given below),
 if for each bounded sequence $\{z_n\}$, $\{T_{n} z_n\}$ is bounded and 
either $(vi)$ $\sum_{n=0}^\infty \norm{ T_{n+N} z_n - T_n z_n} < \infty$ or 
$(vi')$ $ \norm{ T_{n+N} z_n - T_n z_n} /\alpha_n \to 0$ then \hypothd{1}{N} is 
satisfied (again using previous remark about mixing between conditions with or without prime). 
In the previous case, \hypothd{1}{N} is thus implied by 
\hypothdp{1}{N} which is stated now~:
\end{rem}

\hypothdp{1}{N}: For a fixed given $N \ge 1$ and a given sequence $\{\alpha_n\}$ which satisfies 
\hypothd{3}{N} a sequence of mappings $\{T_n\}$ will be said to verify \hypothdp{1}{N}, if 
given bounded sequence $\{z_n\}$, we have $\norm{ T_{n+N} z_n - T_n z_n} \le \rho_n $ 
with either $(vi)$ $\sum_{n=0}^\infty \rho_n < \infty$ or $(vi')$ $ \rho_n /\alpha_n \to 0$. 

\hypothd{2}{p}: 
For a given $p\in X$, a sequence $\{x_n\}$ will be said to verify \hypothd{2}{p} if we have 
\begin{equation}
  \limsup_{n\to \infty} \psca{ f(p) -p ,  J(x_{n} -p) } \le 0\,.
\end{equation}

\hypothd{3}{N}: 
For a fixed given $N \ge 1$, a  sequence of real numbers $\{\alpha_n\}$ will be said to verify
\hypothd{3}{N} if the sequence $\{\alpha_n\}$ is such that $(i)$
$\alpha_n \in (0,1)$, $(ii)$ $\lim_{n\to \infty} \alpha_n = 0$, $(iii)$ 
$\sum_{n=0}^\infty \alpha_n = \infty$ and either $(iv)$ 
$\sum_{n=0}^\infty |\alpha_{n+N}-\alpha_n| < \infty$ 
or $(iv')$ $\lim_{n\to \infty} (\alpha_{n+N}/\alpha_n) = 1$.

We can now formulate the main theorem of the paper~: 

\begin{thm}\label{mainresult} Let $X$ be a ${\cal B}$ real Banach space, 
$C$ a closed convex subset of $X$, $T_n : C \mapsto C$ a sequence of nonexpansive 
mapping, $T$ a nonexpansive mapping and $f \in \Pi_C$. We assume that $\Fix(T) \ne \emptyset$ 
and that for all $n\in \N$ $\Fix(T) \subset \Fix(T_n)$. 
Let $\{\alpha_n\}$ be a sequence of real numbers for which there exists a fixed $N\ge 1$ 
such that \hypothd{3}{N} is satisfied and suppose that there exists $p\in \Fix(T)$ such 
that \hypothd{2}{p} is satisfied, then the sequence $\{x_n\}$ defined by (\ref{seqx}) 
converges strongly to $p$. 
\end{thm}

\begin{myproof} The proof uses a set of Lemmas which are given in section \ref{secbib}.
Since $p$ is in $\Fix(T_n)$ for all $n$ we can use Lemma \ref{lemma1} 
to obtain the boundedness of the sequence $\{x_n\}$. Thus we can conclude using 
Lemma \ref{conclusion}. 
\end{myproof}

\begin{cor}\label{cormainresult}Assume that the hypothesis of Theorem \ref{mainresult} except 
\hypothd{2}{p} are satisfied. Suppose that \hypothd{1}{N} or \hypothdp{1}{N} is satisfied and that for 
each bounded sequence $\{y_n\}$, the sequence $\norm{T_n y_n - T y_n } \rightarrow 0$. 
Then the conclusion of Theorem \ref{mainresult} remains for $p={\mathbf Q}(f)$.
\end{cor}

\begin{myproof} We just need to prove that \hypothd{2}{p} is satisfied for $p={\mathbf Q}(f)$. We first 
show that if \hypothdp{1}{N} is satisfied then \hypothd{1}{N} is also satisfied. As in previous theorem 
$\{x_n\}$ is a bounded sequence. Then, let $p \in \Fix(T)$, we have~: 
\begin{eqnarray}
  \norm{T_n x_n -Tx_n} 
  &\le&  \norm{T_n x_n - T_n p } +  \norm{T_n p - T p } +  \norm{Tp - T x_n } \nonumber \\
  &\le&  2\norm{x_n -  p} + \norm{T_n p - T p }. \nonumber
\end{eqnarray}
and since $\norm{T_n p -T p} \rightarrow 0$ by hypothesis we have that $\{T_n(x_n)\}$ is bounded. 
As shown in remark \ref{hypoth1rem}  we are within the case where \hypothd{1}{N} is implied 
by \hypothdp{1}{N}. Applying Lemma \ref{lemma2} and Corollary 
\ref{coro1} we obtain the convergence of $\norm{Tx_n -x_n}$. 
We can then apply Lemma \ref{limsuplemma} to obtain \hypothd{2}{p} for $p={\mathbf Q}(f)$. 
\end{myproof}

Corollary \ref{cormainresult} can be extended as follows when a constant $T$ cannot be found. 

\begin{cor}\label{cormainresultplus}Assume that the hypothesis of Theorem \ref{mainresult} except 
\hypothd{2}{p} are satisfied. Suppose that \hypothdp{1}{N} is satisfied and that $\{T_n x_n\}$ is bounded and 
that from each subsequence $\sigma{(n)}$ we can extract   a subsequence $\mu{(n)}$ and find a fixed 
mapping $T_\mu$ such that $$\norm{T_{\mu(n)} x_{\mu(n)} - T_\mu x_{\mu(n)}} \to 0.$$
If $F=\Fix{(T_\mu)}$ does not depend on $\mu$, then the conclusion of 
Theorem \ref{mainresult} remains for $p={\mathbf Q}_F(f)$. 
\end{cor}

\begin{myproof} We just need to prove that \hypothd{2}{p} is satisfied for $p={\mathbf Q}(f)$. Using remark 
  \ref{hypoth1rem} we are in the case where \hypothd{1}{N} is implied by \hypothdp{1}{N}.
  Using \hypothd{1}{N} we first easily obtain that $\norm{x_n -T_n x_n}\to 0$ by an argument 
  similar to Corollary \ref{coro1}. Then \hypothd{2}{p} for $p={\mathbf Q}(f)$ follows from Corollary \ref{extension}.
\end{myproof}

We can now consider the case of composition. Assume that $\{T_n^1\}$ and $\{T_n^2\}$ 
satisfy \hypothdp{1}{N} with sequences denoted by $\rho^i_n$. Assume also 
that for a bounded sequence $\{z_n\}$ then the sequences $\{ T_{n+N}^2 z_n\}$ and 
$\{T_{n+N}^1 T_{n+N}^2 z_n\}$ and also bounded. Then  it is straightforward, since the mappings 
$T^1_n$ are nonexpansive, that~: 
\begin{eqnarray}
\norm{ T_{n+N}^1 T_{n+N}^2 z_n - T_{n}^1 \cdots T_{n}^2 z_n} 
  & \le & \rho^1_n + \norm{ T_{n+N}^2 z_n -  T_{n}^2 z_n}\,.   \nonumber
\end{eqnarray}
Thus the composition $T^1_n \circ T^2_n$ satisfy \hypothdp{1}{N} with $\rho_n \defpar \rho^1_n+\rho^2_n$. 
This lead us to propose the following Corollary for dealing with composition~:

\begin{cor}\label{corcomp}Assume that the hypothesis of corollary \ref{cormainresultplus} are 
satisfied  for the sequence $\{T_n^1\}$ with \hypothdp{1}{N} and for $\{T_n^2\}$ also with \hypothdp{1}{N}. 
Then the conclusion of Theorem \ref{mainresult} remains for the sequence $\{T_n^1\circ T_n^2\}$ 
with  $p={\mathbf Q}_F(f)$ and $F=Fix(T^1_\mu \circ T^2_\rho)$. 
\end{cor}

\begin{myproof}As pointed out before the statement of the corollary the composition $T^1_n \circ T^2_n$ 
satisfy \hypothdp{1}{N}. Consider a subsequence $\sigma{(n)}$ we can find first a subsequence 
$\mu_2(n)$ and $\mu_2$ such that~: 
$$\norm{T^2_{\mu(n)} x_{\mu(n)} - T^2_\mu x_{\mu(n)}} \to 0.$$
Then, using properties of the $T^1_n$ sequence, we can re-extract a new subsequence $\rho(n)$ and 
$\rho$ such that~: 
$$\norm{T^1_{\rho(n)} T^2_{\rho(n)} x_{\rho(n)} - T^1_\rho T^2_{\rho(n)} x_{\rho(n)}} \to 0.$$
Since we have~: 
\begin{eqnarray}
  \norm{T^1_{\rho(n)} T^2_{\rho(n)} x_{\rho(n)} - T^1_\rho T^2_\mu  x_{\rho(n)}} 
  & \le & \norm{T^1_{\rho(n)} T^2_{\rho(n)} x_{\rho(n)} - T^1_\rho T^2_{\rho(n)} x_{\rho(n)}} \nonumber \\ 
  &+&  \norm{T^2_{\rho(n)} x_{\rho(n)} - T^2_\rho x_{\rho(n)}} \nonumber 
\end{eqnarray}
When obtain the conclusion for the composition. 
\end{myproof}

Recall that a mapping $T$ is {\em attracting non expansive} if it is nonexpansive and satisfies~: 
\begin{eqnarray}
  \norm{T x - p } < \norm{x -p} \mbox{ for all } x \not\in \Fix{T} \text{ and } p \in \Fix{T}.
\end{eqnarray}
In particular a {\em firmly nonexpansive} mapping, {\em i.e} $\norm{T x - Ty }^2 \le \psca{x-y,T x - Ty}$ 
is attracting nonexpansive \cite{goebel-reich}. 

\begin{rem} In the previous corollary, we obtain a fixed point of a composition and 
in practice the aim is to obtain a common fixed point of two mappings. 
If the mappings $T^1_\mu$ and $T^2_\rho$ are attracting, have a common fixed point 
and $T^1_\mu$ or $T^2_\rho$ is attracting then we will have 
$\Fix{T^1_\mu}\cap \Fix{T^2_\rho} = \Fix{T^1_\mu \circ T^2_\rho}$. The proof is contained 
in \cite[Proposition 2.10 (i)]{ba-bo} and given in Lemma \ref{lemmcompose} for completeness. 
\label{remcompose}
\end{rem}

\begin{rem}
Note that if $X$ is a strictly convex Banach space, then  for $\lambda\in(0,1)$ the mapping 
$T_\lambda \defpar (1-\lambda)I + \lambda T$ is attracting nonexpansive when $T$ is nonexpansive. 
Extension to a set of $N$ operators is immediate by induction. This gives a way to build 
attracting nonexpansive mappings and mixed with previous remark it gives \cite[Proposition 3.1]{takahashi}.
\end{rem}

\begin{rem}Note also that, when $X$ is strictly convex, an other way to 
obtain $F= \cap_i \Fix{(T_i)}$ for a sequence of nonexpansive mappings $\{T_i\}$ is to use 
$T= \sum_{i} \lambda_i T_i$  with a sequence $\{\lambda_i\}$ of real positive numbers such 
that $\sum_{i} \lambda_i =1$ \cite[Lemma 3]{bruck}.
\label{rembruck}
\end{rem}

\subsection{Example $1$} 

\begin{thm} \cite[Theorem 4.2]{xu-2004}
 Let $X$ be a {\cal B} real Banach space, $C$ a closed convex subset of $X$,
 $T :C\mapsto C$ a nonexpansive mapping with $\Fix(T )\ne \emptyset$, and 
 $f$ an $\alpha$-contraction. Then when the sequence $\{\alpha_n\}$ satisfies 
 \hypothd{3}{1} the sequence $\{x_n\}$ 
 defined by \eqref{seqx} with $T_n \defpar T$ converges strongly to ${\mathbf Q}(f)$. 
 \label{thmun}
\end{thm} 

\begin{myproof} Here the sequence $T_n$ does not depend on $n$. We just apply Corollary \ref{cormainresult} 
to get the result. Of course, if the sequence $\{x_n\}$ is bounded then $\{T_n(x_n)=Tx_n\}$ is bounded 
and equation \eqref{hunequation} of \hypothd{1}{1} is then satisfied with 
$\delta_n = |\alpha_n - \alpha_{n+1}|$. Since $\{\alpha_n\}$ satisfies \hypothd{3}{1}, $\{\delta_n\}$ 
satisfies \hypothd{1}{1}. 
We also have $\norm{T_n x_n -T x_n}= 0 \to 0$ and the conclusion follows. 
\end{myproof}

\begin{rem}Suppose now that $T \defpar \sum_{i} \lambda_i T_i$ where $\{\lambda_i\}$ is a sequence of 
positive real numbers such that $\sum_{i} \lambda_i=1$ and the $T_i$ mappings are all supposed 
nonexpansive. Then, we can apply Theorem \ref{thmun} to obtain the strong convergence  of 
the sequence $\{x_n\}$ to ${\mathbf Q}_{\Fix{T}}(f)$. Moreover, If we assume that $X$ is strictly convex 
then using remark \ref{rembruck} we obtain a strong convergence to ${\mathbf Q}_{F}(f)$ with 
$F \defpar \cap_{i\in I} \Fix(T_i)$. 
\end{rem}

This can be extended to the case when the $\lambda_i$ also depends on $n$ and recover \cite[Theorem 4]{kimura} as follows~:

\begin{cor} Let $X$ be a strictly convex {\cal B} real Banach space, $C$ a closed convex subset of $X$,
 $T_i :C\mapsto C$ for $i \in I$ a finite family of nonexpansive mapping with $\cap_{i\in I} \Fix(T_i)\ne \emptyset$, and 
 $f$ an $\alpha$-contraction. For a sequence $\{\alpha_n\}$ satisfying  \hypothd{3}{1} we consider the 
sequence $\{x_n\}$  defined by \eqref{seqx} with $T_n \defpar \sum_{i\in I} \lambda_{i,n} T_i$. 
Assume that for all $i$ and $n$ $\lambda_{i,n} \in [a,b]$ with $a >0$ and $b < \infty$ either 
$\sum_n \lambda_{i,n} < \infty$ or $\lambda_{i,n}/\alpha_n \to 0$ then $\{x_n\}$ converges strongly to ${\mathbf Q}_F(f)$ 
with $F= \cap_{i \in I} \Fix(T_i)$ 
\label{corfinitesum}
\end{cor}

\begin{myproof}The proof is given by an application of corollary \ref{cormainresultplus}. 
Indeed since the $\lambda_{i,n}$ are bounded $T_n x_n$ remains bounded for a bounded sequence $x_n$. 
Then $T_n$ satisfies \hypothdp{1}{1} with $\rho_n = \sum_{i\in I} \lambda_{i,n}$. 
By extracting from each given subsequence $\sigma{(n)}$ a subsequence $\mu{(n)}$ such that 
$\lim_{n\to \infty } \lambda_{i,\mu(n)} = \overline{\lambda}_i $ for all $i\in I$ we can use corollary \ref{cormainresultplus}. 
Finally, noting that, for a strictly convex space $X$, the fixed points of 
$T_{\overline{\lambda}} \defpar \sum_{i\in I} \overline{\lambda}_i T_i$ does not depend on $\overline{\lambda}$ and is equal 
to $\cap_{i\in I} \Fix(T_i)$ we conclude the proof.
\end{myproof}

\subsection{Example $1'$} 

In \cite{song-chen-1} The following algorithm is considered~: 
\begin{equation}
  y_{n+1} = P (\alpha_n f(y_n) + (1-\alpha_n) T y_n)
  \label{seqxsc1}
\end{equation}
Where $P : X \mapsto C$ is a sunny nonexpansive retraction, $f : C \mapsto X $ an $\alpha$-contraction and 
$T : C \mapsto X $ a nonexpansive mapping such that $\Fix(T) \ne \emptyset$. 

If we consider the sequence $x_{n+1} = \alpha_n f(y_n) + (1-\alpha_n) T y_n$ then we have $y_{n+1} = Px_{n+1}$ and thus 
\begin{equation}
  x_{n+1} = \alpha_n f(P(x_n)) + (1-\alpha_n) T(P(x_n))
  \label{seqxsc1p}
\end{equation}
Since $f\circ P$ is an $\alpha$-contraction from $X$ onto $X$ and $T\circ P$ a non expansive mapping 
from $X$ onto $X$ we can use the previous theorem to obtain the strong convergence of the sequence $\{x_n\}$ 
to $x$ a fixed point of $T\circ P$ such that $x= P_{Fix{(T\circ P)}} f(T(x))$ ($P_S$ is the metric projection on $S$). 
We thus obtain now the strong convergence of the initial sequence $\{y_n\}$ to $y=P(x)$ and since 
$x$ is a fixed point of $T\circ P$, $y$ is a fixed point of $P \circ T$. 

If we suppose in addition that $X$ is such that $J$ (or $J_\phi$) is norm-to-weak$^\star$ continuous 
(i.e $X$ is smooth) and that $T$ satisfy the {\rm weakly inward condition} then we can use 
the result of \cite[Lemma 1.2]{song-chen-1} which state that $Fix(T) = Fix(P \circ T)$ to conclude that 
$y$ is in fact a fixed point of $T$ and recover the result of \cite[Theorem 2.4]{song-chen-1}. 

\subsection{Example $2$} 
We consider now the example given in \cite{xu-2005} where the sequence $\{x_n\}$ is given by~: 
\begin{eqnarray} 
  y_n = \beta_n x_n + (1-\beta_n) Tx_n \label{Mann} \nonumber \\ 
  x_{n+1} = \alpha_n u + (1- \alpha_n) y_n \nonumber
\end{eqnarray}
With a sequence of mappings $T_n x \defpar \beta_n x + (1-\beta_n) T x$. This problem is 
rewritten as follows~: 
\begin{equation}
  x_{n+1} = \alpha_n f(x_n) + (1-\alpha_n) T_n x_n 
  \label{seqxman}
\end{equation}

\begin{thm} Let $X$ be a ${\cal B}$ real Banach space, $C$ a closed convex subset of $X$,
 $T :C\mapsto C$ a nonexpansive mapping with $\Fix(T )\ne \emptyset$, and 
 $f$ an $\alpha$-contraction. When the sequence $\{\alpha_n\}$ satisfies 
 \hypothd{3}{1} and the sequence $\{\beta_n\}$ converges to zero and satisfy either 
 $\sum_{n=0}^\infty | \beta_{n+1} - \beta_{n} | < \infty$ or $| \beta_{n+1} - \beta_{n} |/\alpha_n \to 0$.
 Then, the sequence $\{x_n\}$  defined by \eqref{seqxman} converges strongly to ${\mathbf Q}(f)$. 
 \label{exampledeux}
\end{thm}

This theorem is very similar to \cite[Theorem 1]{xu-2005} where $f$ was supposed to be constant. It could 
be covered by corollary \ref{corfinitesum} but here strict convexity is not needed. 

\begin{myproof} We easily check that the fixed points $p$ of $T$ are fixed points of $T_n$ for all $n\in \N$ and $T_n$ is 
nonexpansive for all $n$. Thus by Lemma \ref{lemma1} the sequence $\{x_n\}$ is bounded . If 
the sequence $\{x_n\}$ is bounded then $ \norm{T_n(x_n)} \le max({\norm{x_n}, \norm{Tx_n}})\}$ 
is bounded too. Since~: 
\begin{equation}
  \norm{T_n y_n -T y_n} \le \beta_n (\norm{y_n} + \norm{T y_n}) 
\end{equation}
we have $\norm{T_n y_n -T y_n} \rightarrow 0$ for each bounded sequence $\{y_n\}$. 
It is easily checked that \hypothd{1}{1} is satisfied 
with $\delta_n = |\alpha_{n+1}-\alpha_n|+ |\beta_{n+1}-\beta_n|$.
The conclusion follows from Corollary \ref{cormainresult}.
\end{myproof}

\subsection{Example $3$} 

We consider here the accretive operators example given in \cite{xu-2005} or \cite{xu-2006}~:
\begin{equation}
  x_{n+1} = \alpha_n f(x_n) + (1-\alpha_n) T_n x_n 
  \label{seqxmanacret}
\end{equation}
Where $T_n x= J_{r_n} x$ and $J_\lambda$ is the resolvent of an $m$-accretive operator $A$, 
$J_\lambda x = ( I + \lambda A)^{-1}$. 
The following theorem is similar to \cite[Theorem 4.2, Theorem 4.4]{xu-2006} or 
\cite[Theorem 2]{xu-2005}.

\begin{thm} Let $X$ be a ${\cal B}$ real Banach space, $A$ an $m$-accretive operator in $X$ 
  such that $A^{-1}(0)\ne \emptyset$. We assume here that $C\defpar \overline{D(A)}$ where $D(A)$ is 
the domain of $A$ and suppose that $C$ is convex. Suppose that \hypothd{3}{1} is satisfied by the sequence 
$\{ \alpha_n \}$ and that the sequence $r_n$ is such that $r_n \ge \epsilon >0$ and either 
$\sum_{0}^\infty |1 - r_n/r_{n+1}| < \infty$ or $|1 - r_n/r_{n+1}|/\alpha_n \to 0$, then the sequence $\{x_n\}$ 
 defined by \eqref{seqxmanacret} converges strongly to a zero of $A$. 
\end{thm}

\begin{myproof} We first note that \cite[p 632]{xu-2006}, for $\lambda > 0$, $\Fix{(J_\lambda)} = F$ 
where $F$ is the set of zero of $A$ and for an $m$-accretive operator $A$, $J_\lambda$ is non expansive from 
$X \mapsto \overline{D(A)}$. Using the resolvent identity 
$J_\lambda x = J_\mu ( (\mu/\lambda) x + (1- \mu/\lambda) J_\lambda x)$ we obtain~: 
\begin{equation}
  \norm{T_{n+1} z_n -T_n z_n} \le \left| 1 - \frac{r_n}{r_{n+1}} \right| ( \norm{z_n} + \norm{T_n z_n}) 
\end{equation}
and since the sequence $T_n y_n$ is bounded for a bounded sequence $y_n$ 
(for $p \in A^{-1}(0)$ we have $\norm{T_n y_n -p }\le \norm{y_n -p}$) we can apply remark 
\ref{hypoth1rem} in order to obtain \hypothd{1,1}. 
We thus have $\norm{x_{n+1} - x_n} \rightarrow 0$ by Lemma \ref{lemma2} and 
$\norm{x_n -T_n x_n} \rightarrow 0$ by~: 
\begin{eqnarray}
  \norm{x_{n}- T_n x_{n}} &\le & \norm{x_{n}- x_{n+1} } + \norm{x_{n+1}- T_n x_{n} } \nonumber \\ 
  & \le & \norm{x_{n}- x_{n+1} } + \alpha_n (\norm{f(x_n)} + \norm{T_n(x_n)}) \nonumber
\end{eqnarray}
Take now $r$ such that 
$0 < r < \epsilon$ and define $T \defpar J_r$ then we have~: 
\begin{equation}
  \norm{T_{n} x_n - T x_n } \le \left| 1 - \frac{r}{r_{n}} \right| \norm{x_n -T_n x_n} 
\end{equation}
We thus obtain that $x_n -T x_n \rightarrow 0$ from~: 
\begin{equation}
  \norm{ x_n - T x_n } \le \norm{x_n -T_n x_n} + \norm{T_{n} x_n - T x_n }
\end{equation}
The conclusion is obtained through Corollary \ref{cormainresult}.
\end{myproof}

\subsection{Example $4$} 

We consider here the example given in \cite{song-chen} 
\begin{eqnarray} 
  x_{n+1} = \alpha_n f(x_n) + (1- \alpha_n) T_n y_n \label{seqxper}
\end{eqnarray}
where $T_n = Q_{n \text{ mod } N}$, where $N\ge 1$ is a fixed integer and the $(Q_l)_{l=0,\ldots,N-1}$  is a family 
of nonexpansive mappings. 

\begin{thm} Let $X$ be a ${\cal B}$ real Banach space, $C$ a closed convex subset of $X$,
  $Q_l :C\mapsto C$ for $l \in \{1,\ldots,N\}$ a family of nonexpansive mappings such that 
 $F \defpar \cap_{l=0}^{N-1} \Fix(Q_l )$ is not empty and 
\begin{equation}
  \cap_{l=0}^{N-1} \Fix(Q_l ) = \Fix({T_{n+N}T_{n+N-1}\cdots T_{n+1}}) \text{ for all } n \in \N
\end{equation}
and $f$ an $\alpha$-contraction. When the sequence $\{\alpha_n\}$ satisfies 
\hypothd{3}{N} then the sequence $\{x_n\}$  defined by \eqref{seqxper} converges 
strongly to ${\mathbf Q}_F(f)$.
\end{thm}

\begin{myproof} By Lemma \ref{lemma1}, since the $T_n$ have a common fixed point, 
the sequence $\{x_n\}$ is bounded. Since the sequence of mappings $T_n$ is periodic, 
the sequence $\{T_n x_n\}$ is bounded and equation \eqref{hunequation} of \hypothd{1}{N} 
is obtained for $\delta_n = | \alpha_n - \alpha_{n+N} |$ using \eqref{hunequationrem}.
Since $\{\alpha_n\}$ satisfies \hypothd{3}{N}, $\{\delta_n\}$ satisfies \hypothd{1}{N}. 
Thus, using Lemma \ref{lemma2} we obtain that $\norm{x_{n+N} -x_n} \to 0$. 
Since $\norm{x_{n+1} -T_n x_n}\le \alpha_n (\norm{f(x_n)} + \norm{T_nx_n})$, we have 
$\norm{x_{n+1} -T_n x_n}\to 0$. We introduce the sequence of mappings 
$A_n^{(N,\alpha)} \defpar T_{n+N-1}\cdots T_{n+\alpha}$ for $\alpha \ne N$ and $A_n^{(N,N)}=Id$.
Using Lemma \ref{lemmex3}, given just after this proof, we conclude that~: 
$\norm{ x_{n+N} - A_n^{(N,0)} x_n }\to 0$. This combined with  $\norm{x_{n+N} -x_n} \to 0$ 
gives $\norm{ x_{n+N} - A_n^{(N,0)} x_n }\to 0$. 
Note now that the mappings $A_n^{(N,0)}$ are in finite number are all nonexpansive and 
share common fixed points by hypothesis. Thus we can prove that \hypothd{2}{p} is satisfied for 
$p={\mathbf Q}_F(f)$. 
Let $p={\mathbf Q}_F(f)$ we suppose that \hypothd{2}{p} is not satisfied, then it possible to extract 
a subsequence of $\{x_{\sigma{(n)}}\}$ such that~: 
\begin{equation}
  \lim_{n\to \infty} \psca{ f(p) -p ,  J(x_{\sigma{(n)}} -p) } \le 0 
\end{equation}
But it is then possible to find $q \in\{0,\ldots,N-1\}$ and an extracted new subsequence 
$\mu(n)$ from $\sigma(n)$ such that $\mu_{(n)} \text{ mod } N = q$. 
We thus have $\norm{ x_{\mu{(n)}} - T x_{\mu{(n)}}} \to 0$, 
with $T \defpar A^{(N,0)}_{q}$ which is now a fixed mapping and $\Fix{(T)}=F$. 
Then \hypothu{2}{p} should be true by Lemma \ref{limsuplemma} and this leads to a contradiction. 
The conclusion follows by \ref{conclusion}.
\end{myproof}

\begin{lem} Let $N\in \N$,  $\alpha \in \{0,\ldots,N\}$ and $A_n^{(N,\alpha)} \defpar T_{n+N-1}\cdots T_{n+\alpha}$ for $\alpha \ne N$ and $A_n^{(N,N)}=Id$. Assume that $\norm{x_{n+1} -T_n x_n} \to 0$ then $\norm{ x_{n+N} - A_n^{(N,0)} x_n }\to 0$. 
\label{lemmex3}
\end{lem} 
\begin{myproof} We have for $\alpha \in \{0,\ldots,N-1\}$ by definition of $A_n^{(N,\alpha)}$ 
and using the fact that $A_n^{(N,\alpha)}$ is nonexpansive~: 
\begin{eqnarray}
  \norm{A_n^{(N,\alpha+1)} x_{n+\alpha+1} - A_n^{(N,\alpha)} x_{n+\alpha}} 
    &= & \norm{ A_n^{(N,\alpha+1)} x_{n+\alpha+1} - A_n^{(N,\alpha+1)} T_{n+\alpha} x_{n+\alpha} } 
      \nonumber \\ 
      & \le & \norm{ x_{n+\alpha+1} -  T_{n+\alpha} x_{n+\alpha} }\nonumber
\end{eqnarray}
Thus~: 
\begin{eqnarray}
  \norm{x_{n+N} - A_n^{(N,0)} x_{n}} 
  &\le & 
  \sum_{\alpha=0}^{N-1} \norm{ x_{n+\alpha+1} -  T_{n+\alpha} x_{n+\alpha} } \nonumber
\end{eqnarray}
and the result follows.
\end{myproof}

\subsection{Example $5$} 

Let $\Gamma^{(j)}_n$ for $j\in \{1,\ldots,m\}$ be a sequence of mappings defined recursively as follows~: 
\begin{equation}
  \Gamma^{(j)}_{n} x  \defpar \beta^{(j)}_n x  + (1-\beta^{(j)}_n) T_j \Gamma^{(j+1)}_n x 
  \text{ and } \Gamma^{(m+1)}_n x = x
  \label{gammaj}
\end{equation} 
where the sequences $\{\beta^{(j)}_n\} \in (0,1)$, and $\{T_j\}$ for $j\in\{1,\ldots,m\}$ are 
nonexpansive mappings. We want to prove here the convergence of the sequence generated by the iterations~: 
\begin{equation}
  x_{n+1} = \alpha_n f(x_n) + (1-\alpha_n) \Gamma^{(1)}_n x_n 
  \label{seqxfour}
\end{equation}

\begin{thm} Let $X$ be a ${\cal B}$ real Banach space, $C$ a closed convex subset of $X$,
  $T_j :C\mapsto C$ for $j \in \{1,\ldots,m\}$ a family of nonexpansive mappings such that 
 $\cap_{l=1}^{m} \Fix(T_j )$ is not empty and $f$ an $\alpha$-contraction. 
When the sequence $\{\alpha_n\}$ satisfies \hypothd{3}{N} and for $j\in \{1,\ldots,m\}$ the sequences
$\{\beta^{(j)}_n\}$ satisfy $\lim_{n\to \infty }\beta^{(j)}_{n}=0$ and
 either $\sum_{n=0}^\infty | \beta^{(j)}_{n+1} - \beta^{(j)}_{n} | < \infty$ or 
$| \beta^{(j)}_{n+1} - \beta^{(j)}_{n} |/\alpha_n \to 0$ then the sequence defined by 
\eqref{seqxfour} converges strongly to ${\mathbf Q}_F(f)$ associated to $F=\Fix{(T_1\cdots T_m)}$.
\end{thm}

\begin{myproof} Note first that by an elementary induction $\Gamma^{(1)}_n$ is a nonexpansive mapping. 
If we assume that $p$ is a common fixed point to the mappings $T_i$ then $p$ is a fixed point of the 
mappings $\Gamma^{(j)}_n$. By Lemma \ref{lemma1} the sequence  $\{x_n\}$ is bounded. Then using Lemma \ref{lemma5} 
, given just after this proof, combined with the boundedness of $\{x_n\}$, \hypothd{1}{1} is valid with 
\begin{equation}
  \delta_n =  \sum_{p=1}^m | \beta^{(p)}_{n+1} - \beta^{(p)}_{n} | + | \alpha_{n+1} - \alpha_{n} |
\end{equation}
Now if we can prove that~
\begin{equation}
  \norm{ \Gamma^{(1)}_n x_n - T_1T_2\cdots T_m x_n} \rightarrow 0 
\end{equation}
the conclusion will be given by Corollary \ref{cormainresult}. 
The last assetion can easily be obtained by induction on  
$\norm{ \Gamma^{(j)}_n x_n - T_ j\cdots T_m x_n}$, since we have~:
\begin{eqnarray}
  \norm{ \Gamma^{(j)}_n x_n - T_j \cdots T_m x_n} 
  &\le&  \beta_n^{(j)} (\norm{x_n}+ \norm{T_j \cdots T_m x_n}) \nonumber \\
  && + (1-\beta_n) \norm{ T_j \Gamma^{(j+1)}_n x_n - T_j \cdots T_m x_n} \nonumber \\
  &\le&  \beta_n^{(j)} (\norm{x_n}+ \norm{T_j \cdots T_m x_n}) + \norm{\Gamma^{(j+1)}_n x_n - T_{j+1}\cdots T_m x_n}
  \,.\nonumber 
\end{eqnarray}
\end{myproof}

\begin{rem}For $m=1$ we obtain the same result as Theorem \ref{exampledeux}.
\end{rem} 

\begin{lem}\label{lemma5} Let $\Gamma^{(j)}_n$ be the sequence of mappings defined by (\ref{gammaj})
Then we have for $j\in \{1,\ldots,m\}$~:
\begin{equation}
  \norm{\Gamma^{(j)}_{n+1} x -\Gamma^{(j)}_n{x}} 
  \le \left\{ \sum_{p=j}^m | \beta^{(p)}_{n+1} - \beta^{(p)}_{n} | \right\} K
\end{equation} 
where $K$ is a constant which depends on the mappings $(T_p)_{p\ge j}$ and $x$. 
\end{lem}
\begin{myproof} 
Note first that~:
\begin{equation}
  \norm{\Gamma^{(j)}_{n} x} \le \norm{x} +\norm{T_j(\Gamma^{(j+1)}_{n} x)} 
\end{equation}
which applied recursively shows that $\norm{\Gamma^{(j)}_{n} x}$ is bounded by a constant 
which depends on the mappings $(T_p)_{p\ge j}$ and $x$ and not on $n$. 
Then, using the definition of $\Gamma^{(j)}_n$ we have~: 
\begin{eqnarray}
  \norm{\Gamma^{(j)}_{n+1} x - \Gamma^{(j)}_{n}}
  &\le& |\beta^{(j)}_{n+1}-\beta^{(j)}_n| ( \norm{x} + \norm{T_j \Gamma^{(j+1)} x}) \nonumber \\ 
  && + \norm{T_{j}\Gamma^{(j+1)}_{n+1}(x) - T_{j}\Gamma^{(j+1)}_n(x)}
\end{eqnarray}
since $T_j$ is nonexpansive mappings~:
\begin{equation*}
  \norm{\Gamma^{(j)}_{n+1} x - \Gamma^{(j)}_{n}}
  \le |\beta^{(j)}_{n+1}-\beta^{(j)}_n| ( \norm{x} + \norm{T_j \Gamma^{(j+1)} x}) 
  + \norm{\Gamma^{(j+1)}_{n+1}(x) - \Gamma^{(j+1)}_n(x)}
\end{equation*} 
by recursion and since the last term $\Gamma^{(m+1)}_{n+1}(x) - \Gamma^{(m+1)}_n(x)=0$ we obtain the result. 
\end{myproof}

Note that Lemma \ref{lemma5} remains valid for the sequence 
\begin{equation}
  \Gamma^{(j)}_{n} x  \defpar \beta^{(j)}_n g(x)  + (1-\beta^{(j)}_n) T_j \Gamma^{(j+1)}_n x 
  \text{ and } \Gamma^{(m+1)}_n x = x 
\end{equation} 
if $g$ is a nonexpansive mapping. 

\subsection{Example $6$}  

We consider here the example given in \cite{chen-zhang-fan}
\begin{eqnarray} 
  x_{n+1} = \alpha_n f(x_n) + (1- \alpha_n) T_n x_n \nonumber
  \label{chen-zhang-fan-seq}
\end{eqnarray}
where $T_nx \defpar P_C (x- \lambda_n A x)$ and $P_C$ is the metric projection from $X$ to $C$. 
The aim is to find a solution of the variational inequality problem which is to find 
$x \in C$ such that $\psca{Ax,y-x} \ge 0$ for all $y \in C$. The set of solution of the 
variational inequality problem is denoted by $\mbox{VI}(C,A)$. 
The operator $A$ is said to be 
{\em $\mu$-inverse-strongly monotone} if 
\begin{eqnarray} 
  \psca{x-y, Ax -Ay } \ge \mu \norm{Ax -Ay}^2 \text{ for all } x,y \in C \nonumber
\end{eqnarray}
The next theorem is similar to \cite[Proposition 3.1]{chen-zhang-fan}.

\begin{thm} Let $X$ be a real Hilbert space, $C$ a nonempty closed convex, 
$f$ an $\alpha$-contraction, and 
let $A$ be a $\mu$-inverse-strongly 
monotone mapping of $H$ into itself such that $\mbox{VI}(C,A) \ne \emptyset$. 
Assume that \hypothd{3}{1} is satisfied 
and that $\{\lambda_n\}$ is chosen so that $\lambda_n \in [a,b]$ for some $a$, $b$ with 
$0 < a < b < 2 \mu$ and $\sum_{n=1}^\infty | \lambda_{n+1} - \lambda_n| < \infty$. 
then the sequence $\{x_n\}$ generated by \eqref{chen-zhang-fan-seq} converges 
strongly to ${\mathbf Q}_F(f)$ associated to $F=\Fix{(T_\lambda)}$ 
where  $T_\lambda(x) \defpar P_C (x- \lambda A x)$. $F=\Fix{(T_\lambda)}$ does not depend on 
$\lambda$ for $\lambda >0$ and equals $\mbox{VI}(C,A)$.
\end{thm}

\begin{myproof}For $\lambda > 0$, let $T_\lambda x \defpar  P_C (x- \lambda A x)$. When 
$X$ is an Hilbert space we have $\Fix(T_\lambda) = \mbox{VI}(C,A)$. When $A$ is $\mu$-inverse-strongly 
monotone then for, $\lambda \le 2 \mu$,  $I -\lambda A$ is nonexpansive. Thus the mappings $T_n$ 
are non expansive and $\Fix{(T_n)}=  \mbox{VI}(C,A) \ne \emptyset$. 
By Lemma \ref{lemma1} the sequence $\{x_n\}$ is bounded. Since 
$\norm{T_n z} \le K ( \norm{z} + 2\mu \norm{Az})$, the sequence $\{T_nx_n\}$ is bounded too. 
We also have $\norm{ T_{n+1} z_{n} - T_n z_{n}}\le | \lambda_{n+1} - \lambda_{n}|  \norm{ A z_{n}}$ which 
gives \hypothd{1}{N} with $\delta_n = | \lambda_{n+1} - \lambda_{n}| + | \alpha_{n+1} - \alpha_{n}|$ by 
remark \ref{hypoth1rem}. The result follows now from Corollary \ref{cormainresultplus}. Indeed, 
since $\lambda_{\sigma(n)} \in [a,b]$ it is possible to extract a converging subsequence 
$\lambda_{\mu(n)} \to \overline{\lambda} \in [a,b]$ and we then have 
$\norm{ T_{\mu{(n)}} z  - T_{\overline{\lambda}} z } \le | \lambda_{\mu{(n)}} - {\overline{\lambda}}| \norm{Az}$.
Thus $\norm{ T_{\mu{(n)}} x_{\mu{(n)}}  - T_{\overline{\lambda}} x_{\mu{(n)}} }\rightarrow 0$. 
\end{myproof} 

\begin{rem}We can note that for $\lambda < 2\alpha$, $I-\lambda A$ is in fact attracting 
nonexpansive since~:
$$ \norm{ (I-\lambda A) x - (I-\lambda A)y } \le 
  \norm{x -y } + \lambda(\lambda -2 \alpha) \norm{Ax -Ay}^2 .$$
Thus it is also the case for $P_C \circ (I-\lambda A)$ \cite{ba-bo}. 
For a nonexpansive mapping $S$ we can consider the previous theorem with 
$T_\lambda x \defpar S \circ P_C (x- \lambda A x)$ and using Remark \ref{remcompose} (an Hilbert space is 
strcitly convex) to obtain a strong convergence to a point in $\Fix{(T_\lambda)}= \Fix{S} \cap \mbox{VI}(C,A)$ 
and thus fully recover \cite[Proposition 3.1]{chen-zhang-fan}
\end{rem}

\subsection{Example $7$} 

We consider here the equilibrium problem for a bifunction $F: C \times C \mapsto \R$ where $C$ is 
a closed convex subset of a real Hilbert space $X$. The problem is to find $x \in C$ such that $F(x,y) \ge 0$ 
for all $y\in C$. The set of solutions if denoted by $\mbox{EP}(F)$. It is proved in \cite{flam} (See also \cite{combettes}) that 
for $r > 0$, the mapping $T_r : X \mapsto C$ defined as follows~: 
\begin{equation}
  T_r(x) \defpar \left\{ z \in C : F(z,y) + \frac{1}{r} \psca{ y-z,z-x} \ge 0, \forall y \in C \right\}
\end{equation}
is such that $T_r$ is singled valued, firmly nonexpansive ({\em i.e} $ \norm{T_rx -T_r y}^2 \le 
\psca{ T_r x -T_r y, x-y}$ for any $x,y \in X$), $\Fix(T_r)= \mbox{EP}(F)$ and $\mbox{EP}(F)$ is closed and 
convex if the bifunction $F$ satisfies $(A_1) F(x,x)=0$ for all $x\in C$, $(A_2) F(x,y)+F(y,x) \le 0$ for all 
$x$, $y \in C$, $(A_3)$ for each $x$, $y$, $z \in C$ $\lim_{t \to 0} F(tz + (1-t)x ,y ) \le F(x,y)$ and 
$(A_4)$ for each $x\in C$ $y\mapsto F(x,y)$ is convex and lower semicontinuous. 

we can now consider the sequence $\{x_n\}$ given by~: 
\begin{eqnarray} 
  x_{n+1} = \alpha_n f(x_n) + (1- \alpha_n) T_n x_n \nonumber
  \label{plub}
\end{eqnarray}
where $T_n \defpar T_{r_n}$ for a given sequence of real numbers $\{r_n\}$. 

\begin{thm} Let $X$ be a real Hilbert space,$C$ a nonempty closed convex, 
$f$ an $\alpha$-contraction, assume that $\mbox{EP}(F)\ne \emptyset$,  \hypothd{3}{1} is satisfied 
and the sequence $\{r_n\}$ is such that $\liminf_{n\to \infty}  r_n >0$ and either 
$\sum_{n} | r_{n+1} - r_n | < \infty$ of $| r_{n+1} - r_n |/\alpha_n \to 0$. Then,
the sequence $\{x_n\}$ generated by \eqref{plub} converges 
strongly to ${\mathbf Q}_{\mbox{EP}(F)}(f)$. 
\end{thm}

\begin{myproof} Since the $r_n$ are strictly positive the mappings $T_{r_n}$ are non expansive and share 
the same fixed points $\mbox{EP}(F)$ which was supposed non empty. 
By Lemma \ref{lemma1} the sequence $\{x_n\}$ is bounded. 

Using the definition of $T_r(x)$ and the monotonicity of $F$ ($A_2$) easy computations leads to the 
following inequality \cite[p 464]{plubtieng}~: 

\begin{equation}
  \norm{T_r(x) - T_s(y)} \le \norm{x-y} + \left| 1 - \frac{s}{r} \right| \norm{T_r(y) -y} 
\end{equation}
Using $r>0$ such that $r_n > r$ for all $n\in N$ and $y\in \Fix{(T_r)}$ we obtain 
$\norm{T_{r_n}(x_n) - T_r(y)} \le \norm{x_n-y}$ which gives the boundedness of the sequence 
$\{T_{r_n}(x_n)\}$. Moreover, for a bounded sequence $\{y_n\}$ we obtain~: 
\begin{equation}
  \norm{T_{r_{n+1}}(y_n) - T_{r_n}(y_n)} \le \frac{| r_{n+1} - r_n |}{r} \norm{T_{r_n}(y_n) -y_n}
\end{equation}
We thus obtain \hypothd{1}{1} with $\delta_n= | r_{n+1} - r_n | + |\alpha_{n+1} -\alpha_n| $ using remark \ref{hypoth1rem}. 
The result follows now from Corollary \ref{cormainresultplus}. Indeed, 
since $r_{\sigma(n)} > r $ it is possible to extract a converging subsequence 
$r_{\mu(n)} \to \overline{r} > r $ and we then have 
$\norm{ T_{r_\mu{(n)}} z  - T_{\overline{r}} z } \le | r_{\mu{(n)}} - {\overline{r}}| K$. 
Thus $$\norm{ T_{r_\mu{(n)}} x_{\mu{(n)}}  - T_{\overline{r}} x_{\mu{(n)}} } \rightarrow 0\,.$$ 
\end{myproof}

\section{A collection of Lemma}
\label{secbib}

The first Lemma can be used to derive boundedness of the sequence $\{x_n\}$ 
generated by \ref{seqx}.

\begin{lem}\label{lemma1}Let $\{x_n\}$, the sequence generated by the iterations 
\begin{equation}
  x_{n+1} = \alpha_n f(x_n) + (1-\alpha_n) T_n x_n 
  \label{seqx}
\end{equation}
where $f$ is contraction of parameter $\alpha$, $T_n$ is a family of nonexpansive 
mappings and  $\alpha_n$ is a sequence in $(0,1)$.
Suppose that there exists $p$ a common fixed point of $T_n$ for all $n\in \N$. Then, 
the sequence $\{x_n\}$ is bounded. 
\end{lem}
\begin{myproof} The proof exactly follows the proof of \cite[theorem 3.2]{xu-2004}, the only 
difference is that here the mappings $T_n$ are indexed by $n$ but it does not change the proof. 
Obviously we have~: 
\begin{eqnarray}\label{def28}
  \norm{x_{n+1}-p} &\le & \alpha_n \norm{ f(x_n) - p } + (1-\alpha_n) \norm{T_n x_n - p}\nonumber  \\ 
  & \le & \alpha_n \left(\alpha \norm{x_n-p} + \norm{f(p)-p} \right)  + (1-\alpha_n) \norm{x_n-p}\nonumber  \\ 
  & \le & ( 1 - \alpha_n(1-\alpha)) \norm{x_n-p} + \alpha_n(1-\alpha)\frac{\norm{f(p)-p}}{(1-\alpha)}\nonumber\\ 
  & \le & \max \left( \norm{x_n-p}, \frac{\norm{f(p)-p}}{(1-\alpha)} \right).\nonumber
\end{eqnarray}
And, by induction, $\{x_n\}$ is bounded. 
 \end{myproof}

The next lemma aims at proving that the sequence $\{x_n\}$ is asymptotically regular 
{\it i.e} for a given $N\ge 1$, we have $\norm{x_{n+N}-x_n} \rightarrow 0$. 

\begin{lem}\label{lemma2}With the same assumptions as in Lemma \ref{lemma1} 
and assuming that there exists 
$N \ge 1$ such that \hypothd{1}{N} and \hypothd{3}{N} are fulfilled then, 
for the sequence $\{x_n\}$ given by iterations (\ref{seqx}), we have $\norm{x_{n+N}-x_n} \rightarrow 0$.
\end{lem} 

\begin{myproof} Using the definition of $\{x_n\}$ we have~: 
\begin{eqnarray}
  x_{n+N+1}-x_{n+1} &= & \alpha_{n+N} ( f(x_{n+N}) - f(x_{n})) + (\alpha_{n+N} - \alpha_{n})f(x_{n})\nonumber \\ 
  & & +  (1-\alpha_{n+N}) ( T_{n+N} x_{n+N} -T_{n+N} x_{n}) \nonumber\\
  & & + \left( (1-\alpha_{n+N})T_{n+N} x_{n} -  (1-\alpha_{n}) T_{n} x_{n}\right) \nonumber.
\end{eqnarray}
By Lemma \ref{lemma1} the sequence $\{x_n\}$ is bounded, we can therefore use \hypothd{1}{N} with  $\{x_n\}$. 
Since $\{f(x_n)\}$ is bounded too, we can find three constants such that~: 
\begin{eqnarray}
  \norm{x_{n+N+1}-x_{n+1}} &\le & \alpha_{n+N} \alpha \norm{x_{n+N} - x_{n}} + |\alpha_{n+N} - \alpha_{n}| K_1 \nonumber \\ 
    &  & +  (1-\alpha_{n+N}) \norm{x_{n+N} - x_{n}} + \delta_n M  \nonumber  \\ 
  & \le & ( 1 - (1-\alpha) \alpha_{n+N}) \norm{x_{n+N} - x_{n}} + (|\alpha_{n+N} - \alpha_{n}| + \delta_n) K_2 \nonumber
\end{eqnarray}
The proof then follows easily using the properties of $\alpha_n$ i.e \hypothd{3}{N} and Lemma~\ref{lemmaconvjpc}. 
\end{myproof}

The next step is to prove that we can find a fixed mapping $T$ such that $\norm{x_{n}- T x_n} \rightarrow 0$. 
The next corollary gives a simple example for which the property can be derived from Lemma \ref{lemma2}. 
Indeed,  we have  seen specific proofs in previous sections on illustrated examples. 

\begin{cor}Using the same hypothesis as in Lemma \ref{lemma2} and assuming 
  that $\{T_nx_n\}$ is bounded and that 
  $\norm{T_n x_n - T x_n } \rightarrow 0$ we also have $\norm{x_{n}- T x_n} \rightarrow 0$. 
\label{coro1}
\end{cor}
\begin{myproof}
\begin{eqnarray}
  \norm{x_{n}- T x_{n}} &\le & \norm{x_{n}- x_{n+1} } + \norm{x_{n+1}- T x_{n} } \nonumber \\ 
  & \le & \norm{x_{n}- x_{n+1} } + \alpha_n K_1 + (1-\alpha_n) \norm{T_n x_n - T x_n } \nonumber
\end{eqnarray}
and the result follows. 
\end{myproof}

The next Lemma gives assumptions to obtain \hypothd{2}{p} for a given $p$.

\begin{lem}Suppose that $X$ is a ${\cal B}$ real Banach space. 
  Let $T$ be a nonexpansive mapping with $\Fix(T)\ne \emptyset$, $f$ an $\alpha$-contraction 
  and $\{x_n\}$ a bounded sequence 
  such that $\norm{T x_n -x_n}\to 0$. Then for $\tilde{x}={\mathbf Q}(f)$ we have~:  
\begin{equation}
  \limsup_{n\to \infty} \psca{ f(\tilde{x}) -\tilde{x} ,  J(x_{n} -\tilde{x}) } \le 0 
\end{equation}
\label{limsuplemma}
\end{lem} 
\begin{myproof}When $X$ is a ${\cal B}_{\mbox{us}}$ or a ${\cal B}_{\text{rug}}$ the key point is
the fact that $J$ is uniformly norm-to-weak$^\star$ continuous on bounded sets. 

The proof of this Lemma can be found in the proof of Theorem \cite[Theorem 4.2]{xu-2004} or 
\cite[Theorem 3.1]{song-chen}. We just summarize the line of the proof here. 
Let $\tilde{x} \defpar \slim_{t\to 0} x_t$ where $x_t$ solves 
$x_t = tf(x_t) + (1-t) Tx_t$, we thus have~: 
\begin{eqnarray}
  \norm{x_{t} -  x_{n}}^2 & \le & (1-t)^2\norm{ Tx_{t} -  x_{n}}^2 + 2t\psca{ f(x_t) -x_n, J(x_t-x_n) } \nonumber \\ 
  &\le &  (1-t)^2(\norm{Tx_{t} - Tx_{n}} + \norm{Tx_n -x_n})^2 \nonumber \\ 
  && + 2t\psca{ f(x_t) - x_t, J(x_t-x_n) } + 2t  \norm{x_{t} -  x_{n}}^2  \nonumber\\ 
  & \le & (1+t^2)\norm{x_{t} - x_{n}}^2 + a_n(t)  \nonumber\\ 
  &&  +  2t\psca{ f(x_t) - x_t, J(x_t-x_n) } \nonumber\\ 
\end{eqnarray}
where $a_n(t) = 2\norm{Tx_n -x_n}\norm{x_{t} -  x_{n}} + \norm{Tx_n -x_n}^2\to 0$ when $n$ tends to 
infinity. Thus~: 
\begin{eqnarray}
  \psca{ f(x_t) - x_t, J(x_n -x_t) } \le \frac{a_n(t)}{2t} + \frac{t}{2} \norm{x_{t} - x_{n}}^2
\end{eqnarray}
and we have~: 
\begin{eqnarray}
  \lim_{t\to 0} \limsup_{n\to \infty} \psca{ f(x_t) - x_t, J(x_n -x_t) } \le 0
  \label{limsupres}
\end{eqnarray}
We consider now a sequence $t_p \to 0$ and $y_p \defpar x_{t_p}$, then we have 
$y_p \to \tilde{x}$ and with $g(x)\defpar (x)-x$ we have 
\begin{eqnarray}
  \psca{ g(\tilde{x}) , J(x_n - \tilde{x}) } 
    & \le & 
    \psca{ g(y_p) , J(x_n - y_p)} \nonumber \\
    &+ & | \psca{ g(\tilde{x})  , J(x_n - \tilde{x}) -  J(x_n - y_p)} | + 
   (1+\alpha) \norm{ \tilde{x} - y_p} \norm{x_n -y_p} \nonumber 
\end{eqnarray}

Since $J$ is uniformly norm-to-weak$^\star$ continuous on bounded sets and  $y_p \to \tilde{x}$, 
for $\epsilon >0$, we can find $\tilde{p}$ such that for all $p \ge \tilde{p}$ and all 
$n\in \N$ we have~: 
\begin{eqnarray}
  \psca{ g(\tilde{x})  , J(x_n - \tilde{x})}
    & \le & 
    \psca{ g (y_p) , J(x_n - y_p)} + \epsilon 
      (1+\alpha) \norm{ \tilde{x} - y_p} \norm{x_n -yp} 
\end{eqnarray}
Thus~: 
\begin{eqnarray}
  \limsup_{n\to \infty}\psca{ g(\tilde{x}) , J(x_n - \tilde{x})}
    & \le & 
    \limsup_{n\to \infty} \psca{ g(y_p)  , J(x_n - y_p)} + \epsilon 
      + \norm{ \tilde{x} - y_p} K \nonumber\\ 
    & \le & 
      \lim_{p\to \infty}( \limsup_{n\to \infty} \psca{ g(y_p) , J(x_n - y_p)} + \epsilon 
      \norm{ \tilde{x} - y_p} K) 
       \le  \epsilon \nonumber
\end{eqnarray}

Suppose now that $X$ is a ${\cal B}_{\text{rwsc}}$. 
We follow the proof of [Theorem 2.2]{song-chen-1} or \cite[Theorem 3.1]{xu-2006}. 
Let $\tilde{x}={\mathbf Q}(f)$ and consider a subsequence $\{x_{\sigma{(n)}}\}$ such that 
$ \limsup_{n\to \infty} \psca{ f(\tilde{x}) -\tilde{x} ,  J(x_{n} -\tilde{x}) } = 
\lim_{n\to \infty} \psca{ f(\tilde{x}) -\tilde{x} ,  J(x_{\sigma{(n)}} -\tilde{x}) }$. 
It is then possible to re-extract a subsequence $x_{\mu{(n)}}$ weakly 
converging to $x^\star$. Since we have $x_{\mu{(n)}} -T x_{\mu{(n)}} \rightarrow 0$ then $x^\star \in \Fix(T)$ 
using the key property that $X$ satisfies Opial's condition \cite[Theorem 1]{gossez-dozo} and the fact that 
$I-T$ is demi-closed at zero \cite[Lemma 2.2]{song-chen}. 
Thus by definition of $\tilde{x}$ we must have $\psca{ f(\tilde{x}) -\tilde{x} ,  J(x^\star -\tilde{x}) }\le 0 $. 
\end{myproof}

\begin{cor}Suppose that $X$ is a ${\cal B}_{\mbox{us}}$, 
  or a ${\cal B}_{\text{rug}}$, or a ${\cal B}_{\text{rwsc}}$.
  let $f$ a contraction and $\{x_n\}$ a bounded sequence such that 
  $x_n -T_n x_n \to 0$. From each subsequence $\sigma{(n)}$ we can extract 
  a subsequence $\mu{(n)}$ and find a fixed mapping $T_\mu$ such that 
  $\norm{T_{\mu(n)} x_{\mu(n)} - T_\mu x_{\mu(n)}} \to 0$. Then, if $F=\Fix{T_\mu}$ does not depend on $\mu$, 
  for $\tilde{x}={\mathbf Q}(f)$ associated to $F$, we have~:  
  \begin{equation}
    \limsup_{n\to \infty} \psca{ f(\tilde{x}) -\tilde{x} ,  J(x_{n} -\tilde{x}) } \le 0 
  \end{equation}
  \label{extension}
\end{cor}

\begin{myproof}The proof is by contradiction using Lemma \ref{limsuplemma}. 
  Assume that the result is false, then we can find a subsequence $\sigma(n)$ such that 
\begin{equation}
  \label{epscond}
  \limsup_{n\to \infty} \psca{ f(\tilde{x}) -\tilde{x} ,  J(x_{\mu{(n)}} -\tilde{x}) } \ge \epsilon > 0
\end{equation} 
by hypothesis we can extract from $\sigma(n)$ a sub-sequence $\mu{(n)}$ such that 
$\norm{T_{\mu(n)} x_{\mu(n)} - T x_{\mu(n)}} \to 0$. Thus, since 
\begin{equation*}
  \norm{x_{\mu{(n)}} -T x_{\mu{(n)}} } \le   \norm{x_{\mu{(n)}} -T_{\mu(n)} x_{\mu{(n)}} } 
  + \norm{T_{\mu(n)} x_{\mu{(n)}} - T x_{\mu{(n)}}},
\end{equation*} 
we have $x_{\mu{(n)}} -T x_{\mu{(n)}}\rightarrow 0$ we can then apply Lemma \ref{limsuplemma} to the sequence 
$\{x_{\mu{(n)}}\}$ and mapping $T_\mu$ to derive that~: 
\begin{equation*}
  \limsup_{n\to \infty} \psca{ f(\tilde{x}) -\tilde{x} ,  J(x_{\mu{(n)}} -\tilde{x}) } \le 0 
\end{equation*}
for $\tilde{x}={\mathbf Q}(f)$ corresponding to $F=\Fix{T_\mu}$ and since $F$ does not depend on $\mu$, this gives 
a contradiction with (\ref{epscond}).
\end{myproof}

The next Lemma helps concluding the proof. 

\begin{lem}\label{conclusion} Assume that the sequence $\{x_n\}$ 
  given by iterations (\ref{seqx}) 
  is bounded and assume that for $p$, a common fixed point of the mappings ${T_n}, 
  $\hypothd{2}{p} is satisfied and that $(i,ii,iii)$ items of \hypothd{3}{N} is also 
  satisfied\footnote{Note that $(i,ii,iii)$ of \hypothd{3}{N} do not use the value of $N$}.
  Then the sequence $\{x_n\}$ converges to $p$.
\end{lem}

\begin{myproof}
\begin{eqnarray}
  \norm{x_{n+1}- p}^2 &\le & ( 1 - \alpha_n)^2 \norm{ T_nx_{n}- p }^2 
  + 2 \alpha_n \psca{ f(x_n) -p,  J(x_{n+1} -p) }  \nonumber \\ 
  &\le & ( 1 - \alpha_n)^2 \norm{ x_{n}- p }^2 + 2 \alpha_n \psca{ f(x_n) - f(p) , J(x_{n+1} -p)}\nonumber \\
  & & + 2 \alpha_n \psca{ f(p) -p ,  J(x_{n+1} -p) }  \nonumber \\ 
  &\le & ( 1 - \alpha_n)^2 \norm{ x_{n}- p }^2 + 2 \alpha_n \alpha \norm{x_n - p}\norm{x_{n+1} -p)}\nonumber \\
  & & + 2 \alpha_n \psca{ f(p) -p ,  J(x_{n+1} -p) }  \nonumber 
\end{eqnarray}
Note that $\norm{x_{n+1} - p} \le \norm{x_{n} - p} + \alpha_n K$ . Thus~: 
\begin{eqnarray}
  \norm{x_{n+1}- p}^2 
  &\le & ( 1 - \alpha_n)^2 \norm{ x_{n}- p }^2 + 2 \alpha_n \alpha \norm{x_n - p}^2 \nonumber \\
  & & + 2\alpha_n^2 K  + 2 \alpha_n \psca{ f(p) -p ,  J(x_{n+1} -p) }  \nonumber \\ 
  &\le & ( 1 - \alpha_n(1-\alpha) + \alpha_n^2) \norm{ x_{n}- p }^2 \nonumber \\
  & & + 2\alpha_n^2 K  + 2 \alpha_n \psca{ f(p) -p ,  J(x_{n+1} -p) }  \nonumber \\ 
\end{eqnarray}
And we conclude with  Lemma~\ref{lemmaconv}.
\end{myproof}

\begin{lem}\label{lemmaconv}.\cite[Lemma 2.1]{xu-2005}
Let $\{s_n\}$ be a sequence of nonnegative real numbers satisfying the property
\begin{equation*}
  s_{n+1} \le (1-\alpha_n) s_n + \alpha_n \beta_n \text{ for } n \ge 0\,,
\end{equation*}
where ${\alpha}_n \in (0,1)$ and ${\beta_n}$ are sequences of real numbers such that~: 
$(i)$ $\lim_{n\to \infty} {\alpha}_n =0$ and $\sum_{n=0}^\infty \alpha_n = \infty$ 
$(ii)$ either $\limsup_{n\to \infty} \beta_n \le 0$ or $\sum_{n=0}^\infty |\alpha_n\beta_n| < \infty$. Then $\{s_n\}$ converges to zero. 
\end{lem}

\begin{cor}\label{lemmaconvjpc} 
Let $\{s_n\}$ be a sequence of nonnegative real numbers satisfying the property
\begin{equation*}
  s_{n+1} \le (1-\alpha_n) s_n + \alpha_n \beta_n + \alpha_n \gamma_n\text{ for } n \ge 0\,,
\end{equation*}
where ${\alpha}_n \in (0,1)$, ${\beta_n}$ and ${\gamma_n}$ are sequences of real numbers such that~: 
$(i)$ $\lim_{n\to \infty} {\alpha}_n =0$ and $\sum_{n=0}^\infty \alpha_n = \infty$ 
$(ii)$ $\limsup_{n\to \infty} \beta_n \le 0$ and $(iv)$ $\sum_{n=0}^\infty |\alpha_n\delta_n| < \infty$. 
Then $\{s_n\}$ converges to zero. 
\end{cor} 

\begin{myproof}The proof is similar to the proof of Lemma \ref{lemmaconv} \cite[Lemma 2.1]{xu-2005}. 
Fix $\epsilon >0$ and $N$ such that $\beta_n \le \epsilon/2$ for $n \ge N$ and 
$\sum_{j=N}^\infty |\alpha_n\delta_n|\le \epsilon/2$ . Then following \cite{xu-2005} we have
for $n > N$~:
\begin{eqnarray}
  s_{n+1} &\le & \prod_{j=N}^n (1-\alpha_j)s_N + \frac{\epsilon}{2}( 1 - \prod_{j=N}^n (1-\alpha_j)) +  \sum_{j=N}^n |\alpha_n\delta_n| \nonumber \\
  &\le & \prod_{j=N}^n (1-\alpha_j)s_N + \frac{\epsilon}{2}( 1 - \prod_{j=N}^n (1-\alpha_j)) + \frac{\epsilon}{2}
\end{eqnarray}
and then by taking the limit sup when $n \to \infty$ we obtain $\limsup_{n\to \infty} s_{n+1} \le \epsilon$. 
\end{myproof}

A contraction is said to be a {\em  Meir-Keeler contraction} (MKC) if for every $\epsilon>0$ there exits $\delta > 0$ 
such that $\norm{x-y} < \epsilon + \delta $  implies $\norm{\Phi(x)-\Phi(y)} < \epsilon$. 

\begin{lem}\cite{suzuki} Suppose that the sequence $\{x_n\}$ defined by equation (\ref{seqx}) strongly converges 
 for an $\alpha$-contraction $f$ (or a constant function $f$) to the fixed point of $P_F\circ f$ then 
 the results remains valid for a Meir-Keeler contraction $\Phi$. 
\label{mkclem}
\end{lem}

\begin{myproof}Suppose that we have proved that (\ref{seqx}) converges for an $\alpha$-contraction $f$ 
to the fixed point of $P_F\circ f$. Then indeed, the result is true when $f$ is a constant mapping. 
Let $\Phi$ be a Meir-Keeler contraction, fix $y\in C$, when $f$ is constant and equal to $\Phi(y)$ then 
$\{x_n\}$ defined by \eqref{seqx} converges to $P_F(\Phi(y))$. If $\Phi$ is a MKC then since $P_F$ is nonexpansive 
$P_F\circ \Phi$ is also MKC (Proposition 3 of \cite{suzuki}) and has a unique fixed point 
\cite{meir-keeler}. We can consider $z=P_F(\Phi(z))$ and consider two sequences~:
\begin{equation}
  x_{n+1} = \alpha_n \Phi(x_n) + (1-\alpha_n) T_n x_n 
\end{equation} 
\begin{equation}
  y_{n+1} = \alpha_n \Phi(z) + (1-\alpha_n) T_n y_n 
\end{equation} 
Of course $\{y_n\}$ converges strongly to $z$. We now prove that $\{x_n\}$ also 
converges strongly to $z$ following \cite{suzuki}. 
Fix $\epsilon >0$, by Proposition 2 of \cite{suzuki}, we can find $r\in (0,1)$ such that 
$\norm{x-y}\le \epsilon$ implies $\norm{\Phi(x)-\Phi(y)} \le r \norm{x-y}$. 
Choose now $N$ such that $\norm{y_n -z} \le \epsilon (1-r)/r$. 
Assume now that for all $n\ge N$ we have $\norm{x_n-y_n} > \epsilon$ then 
\begin{eqnarray}
  \norm{x_{n+1} - y_{n+1}} &\le&
  (1-\alpha_n) \norm{x_n - y_n} + \alpha_n \norm{\Phi(x_n)-\Phi(y_n)} + \alpha_n \norm{\Phi(y_n)-z} \nonumber \\
  &\le & (1 -\alpha_n(1-r)) \norm{x_n - y_n} + \alpha_n \epsilon \nonumber 
\end{eqnarray} 
We cannot use here directly Lemma \ref{lemmaconv} but following the proof of this Lemma we obtain 
that $\limsup \norm{x_n - y_n} \le \epsilon$. 
Assume now that for a given value of $n$ we have $\norm{x_n-y_n}\le \epsilon$. 
Since $\Phi$ is a MKC we have 
$\norm{\Phi(x)-\Phi(y)} \le \max(r\norm{x-y},\epsilon)$ and since we have 
\begin{equation}
  r\norm{x_n -z} \le r\norm{x_n -y_n}  + r \norm{y_n -z} \le \epsilon 
\end{equation}
we obtain 
\begin{equation}
  \norm{x_{n+1} - y_{n+1}} \le (1-\alpha_n) \norm{T_n x_n -T_n y_n} + \alpha_n \max(r\norm{x_n -z},\epsilon) \le \epsilon\,.
\end{equation}
Thus we have in both cases $\limsup_{n\to \infty} \norm{x_n -y_n} \le \epsilon$ and the conclusion follows. 
\end{myproof}

\begin{lem}\cite[Proposition 2.10 (i)]{ba-bo} Suppose that $X$ is strictly convex, $T_1$ 
an attracting non expansive mapping and $T_2$ a non expansive mapping which have a 
common fixed point. Then~:
$$\Fix(T_1 \circ T_2) = \Fix(T_2 \circ T_1) = \Fix(T_2) \cap \Fix( T_1)\,.$$ 
\label{lemmcompose}
\end{lem}

\begin{myproof}We have $ \Fix(T_2) \cap \Fix( T_1) \subset \Fix(T_2 \circ T_1)$ and $\Fix(T_2) \cap \Fix( T_1) \subset \Fix(T_1 \circ T_2)$. Let $x$ be a common fixed point of $T_1$ and $T_2$. 
If $y$, a fixed point of $T_1 \circ T_2$, is such that $y\not \in \Fix(T_2)$ then since 
$T_1$ is attracting non expansive we have~: 
$$
  \norm{y-x } = \norm{ T_1 \circ T_2 (y) -x } < \norm{ T_2 (y) -x} \le \norm{y-x} 
$$ 
which gives a contradiction. Thus $y$ is a fixed point of $T_2$ and then also of $T_1$. 
If now $y$ a fixed point of $T_2 \circ T_1$ and assume that  $y \not \in \Fix(T_1)$ then we have 
$$
 \norm{y-x } = \norm{ T_2 \circ T_1 (y) -x } \le  \norm{ T_1 (y) -x} < \norm{y-x} 
$$ 
which gives also a contradiction and same conclusion. 
\end{myproof}


\end{document}